\documentclass{article}
\usepackage{amsmath, amssymb}
\usepackage{amssymb,pb-diagram,pst-all,graphicx}
\usepackage{amscd}

\usepackage{fancybox}
\newcommand{\ZZ}{\mathbb Z}
\newcommand{\PP}{\mathbb P}
\newcommand{\QQ}{\mathbb Q}
\newcommand{\CC}{\mathbb C}

\newcommand{\mcB}{\mathcal B}
\newcommand{\mcC}{\mathcal C}

\newcommand{\mcQ}{\mathcal Q}
\newcommand{\mcT}{\mathcal T}

\newcommand{\Sub}{\underline {\mathrm {Sub}}}

\newcommand{\MW}{\mathop {\rm MW}\nolimits}

\newcommand{\NS}{\mathop {\rm NS}\nolimits}
\newcommand{\Div}{\mathop {\rm Div}\nolimits}

\newcommand{\Red}{\mathop {\rm Red}\nolimits}

\newcommand{\Sing}{\mathop {\rm Sing}\nolimits}

\newtheorem{thm}{Theorem}[section]
\newtheorem{cor}{Corollary}[section]
\newtheorem{prop}{Proposition}[section]
\newtheorem{lem}{Lemma}[section]
\newtheorem{defin}{Definition}[section]
\newtheorem{exmple}{Example}[section]
\newtheorem{rem}{Remark}[section]
\newtheorem{qz}{Question}[section]

\newcommand{\I}{\mathop {\rm I}\nolimits}

\newcommand{\III}{\mathop {\rm III}\nolimits}

\newcommand{\qed}{\hfill $\Box$}

\newcommand{\proof}{\noindent{\textsl {Proof}.}\hskip 3pt}
\newcommand{\proofend}{\qed \par\smallskip\noindent}

\renewcommand{\thesubparagraph}{\theparagraph.\@arabic\c@subparagraph}
\begin{document}
  
  \begin{center}
  
 {\bf  \Large 
Geometry of bisections of elliptic surfaces \\ and  \\
Zariski $N$-plets II}
\bigskip

\bigskip
\large 
Shinzo BANNAI,
and
Hiro-o TOKUNAGA

%\footnote{Research partly supported by the research grant 22540052
%from JSPS}

\end{center}
\normalsize

\begin{abstract}
In this article, we continue to study the geometry of bisections of certain rational elliptic surfaces. As an application,
we give examples of Zariski $N+1$-plets of degree $2N+4$ whose irreducible components
are an irreducible quartic curve with $2$ nodes and $N$ smooth conics. Furthermore,
by considering the case of $N = 2$ and combining with known results, a new 
Zariski $5$-plet for reduced plane curves of degree $8$ is given. 
\end{abstract}

\section*{Introduction}
 In this article, we continue to study the topology of reducible plane curves via the 
 geometry of elliptic surface as in \cite{bannai16, bannai-tokunaga, tokunaga10, tokunaga14}. As for terminologies about Zariski pairs or Zariski $N$-plets,
 we refer to \cite{act}.
 
 Let $\mcB$ be a reduced plane curve with irreducible decomposition $\mcB = \mcB_1 + \cdots + \mcB_r$. A smooth conic $\mcC$ is said to be a contact conic to $\mcB$ if (i)
 no singular point of $\mcB$ is contained in $\mcB \cap \mcC$ and (ii) for
 $\forall x \in \mcB\cap \mcC$, the intersection multiplicity $I_x(\mcB, \mcC)$ at $x$ is 
 even.
 
 For a smooth conic $\mcC$, let $f_{\mcC}: Z_{\mcC} \to \PP^2$ be the double
 cover branched along $\mcC$. If $\mcC$ is a contact conic to $\mcB$, for
 any irreducible component $\mcB_i$, $f_{\mcC}^*\mcB_i$ is either irreducible or of the
 form $\mcB_i^+ + \mcB_i^-$, where $\mcB_i^- = \sigma_{f_{\mcC}}^*\mcB_i^+$, $\sigma_{f_{\mcC}}$ being the covering transformation of $f_{\mcC}$. For the latter case, we say
 $\mcB_i$ to be splitting or a splitting curve with respect to $f_{\mcC}$. As we have
 already seen in \cite{bannai16, bannai-shirane16, shirane16, tokunaga10},
 whether $\mcB_i$ is splitting or not is a subtle and interesting question
 in the study of the topology of $\mcC + \mcB$. In fact, we have the 
 following:
 
 \begin{prop}\label{prop:split-1}{(\cite[Proposition~1.3]{shirane16}
 Let $\mcB$ be an irreducible plane curve. Let $\mcC_i$ $(i = 1, 2)$ be
 contact conics to $\mcB$. If 
 $f_{\mcC_1}^*\mcB$ is splitting and $f_{\mcC_2}^*\mcB$ is irreducible, then there 
 exists no homeomorphism $h : (\PP^2, \mcC_1 + \mcB) \to
 (\PP^2, \mcC_2 + \mcB)$ with $h(\mcB) = \mcB$.
 }
 \end{prop}

 In \cite{tokunaga10},  when $\mcB$ is an irreducible plane quartic, a criterion
 for $\mcB$ to be splitting with respect to $f_{\mcC}$ or not is given.
 We here recall it as follows:

 Let $\mcQ$ be a  reduced quartic which is not a union of four concurrent lines
 and choose
 a general point $z_o$ of $\mcQ$. We can associate a rational 
elliptic surface $S_{\mcQ, z_o}$ (see
 \cite[2.2.2]{bannai-tokunaga}, \cite[Section 4]{tokunaga14}) to $\mcQ$ and $z_o$,
 which is given as follows:
 
 \begin{enumerate}
 
 \item[(i)] Let $f'_{\mcQ} : S'_{\mcQ} \to \PP^2$ be a double cover branched along $\mcQ$.
 
 \item[(ii)] Let $\mu: S_{\mcQ} \to S'_{\mcQ}$ be the minimal resolution of $S'_{\mcQ}$.
 
 \item[(iii)] The pencil of curves passing through $z_o$ on $\PP^2$ gives rise to a
  pencil $\Lambda_{z_o}$  of curves of genus $1$ with a unique base 
  point $(f'_{\mcQ}\circ\mu)^{-1}(z_o)$.
  
  \item[(iv)] Let $\nu_{z_o}: S_{\mcQ, z_o} \to S_{\mcQ}$ be the resolution of the 
  indeterminancy for the rational map induced by $\Lambda_{z_o}$. We denote the induced morphism 
  $\varphi_{\mcQ, z_o} : S_{\mcQ, z_o} \to \PP^1$, which gives a minimal elliptic fibration.
  The map $\nu_{z_o}$ is a composition of two blowing-ups and the exceptional curve for
  the second blowing-up gives rise to a section $O$ of $\varphi_{\mcQ, z_o}$.
  Note that we have the following diagram:
  \[
\begin{CD}
S'_{\mcQ} @<{\mu}<< S_{\mcQ} @<{\nu_{z_o}}<<S_{\mcQ, z_o} \\
@V{f'_{\mcQ}}VV                 @VV{f_{\mcQ}}V         @VV{f_{\mcQ, z_o}}V \\
\PP^2@<<{q}< \widehat{\PP^2} @<<{q_{z_o}}< (\widehat{\PP^2})_{z_o},
\end{CD}
\]
where $f_{\mcQ, z_o}$ is a double cover induced by the quotient under the involution $[-1]_{\varphi_{\mcQ, z_o}}$ on $S_{\mcQ,z_o}$, which is given by the inversion with respect to the group law 
 on the generic fiber. Furthermore $q$ is a composition of a finite 
 number of blowing-ups so that
 the branch locus becomes smooth and $q_{z_o}$ is a composition of 
 two blowing-ups.
  
 \end{enumerate}

In the following, we always assume that 

\begin{center}
a reduced quartic $\mcQ$ is 
not a union of four concurrent lines.
\end{center}

By our assumption that $\mcC$ is a contact conic, and since $\mcC$ is rational, $(q\circ f_{\mcQ})^*\mcC = \mcC^+ + \mcC^-$. Furthermore, if we assume $z_o \in \mcC$,
$\mcC^{\pm}$ give rise to sections $s_{\mcC^{\pm}, z_o}$ of 
$\varphi_{\mcQ, z_o}$, respectively. Now  we have

\begin{thm}\label{thm:q-residue}{\cite[Theorem~2.1]{tokunaga10} In the above setting, $f_{\mcC}^*\mcQ$ is
 splitting if and only if $s_{\mcC^+, z_o}$ is $2$-divisible in $\MW(S_{\mcQ, z_o})$
 }
 \end{thm}
 
 Note that $s_{\mcC^-, z_o}$ is $2$-divisible if and only if $s_{\mcC^+, z_o}$ is also two divisible,
 as $s_{\mcC^-, z_o} = [-1]s_{\mcC^+, z_o}$. Also our statement in Theorem~\ref{thm:q-residue}
  involves the choice of $z_o$, although the splitting property does not seem to depend on $z_o$. Later, we will resolve this discordance.
  
 As we have seen in  \cite{tokunaga10},  when $\mcQ$ is irreducible and has two nodes 
 ($2$-nodal quartic) or 
 one tacnode and a general point $z_o$, 
 there exist contact conics $\mcC_1,\, \mcC_2$ through $z_o$
 such that
 (i) $f^*_{\mcC_1}\mcQ$ is splitting and
 (ii) $f^*_{\mcC_2}\mcQ$ is irreducible.

 Since any homeomorphism $h : (\PP^2,  \mcC_1 + \mcQ) \to (\PP^2, \mcC_2+
 \mcQ)$ satisfies $h(\mcQ) = \mcQ$, by Proposition~\ref{prop:split-1}, we have:
 
 \begin{prop}\label{prop:zpair-1}{Let $\mcQ$ be a $2$-nodal quartic or
 an irreducible quartic with one tacnode. Let $z_o$ be a general point and 
 let $\mcC_i$ $(i = 1, 2)$ be
 as above. Then
 $(\PP^2, \mcC_1 + \mcQ)$ is not homeomorphic to 
 $(\PP^2,  \mcC_2 + \mcQ)$. In particular, if the combinatorial type
 of $\mcC_i + \mcQ$ $(i = 1, 2)$ are the same, then
 $(\mcC_1 + \mcQ, \mcC_2 + \mcQ)$ is a Zariski pair.
 }
 \end{prop}
 
 \begin{rem}{\rm In \cite{tokunaga10}, we make use of existence/non-existence
 of dihedral covers of $\PP^2$ branched along $\mcC_i + \mcQ$ ($i = 1, 2$) to prove Proposition \ref{prop:zpair-1}.
 }
 \end{rem}
 
 In this article, we first consider a generalization of Proposition~\ref{prop:zpair-1}
 as follows:
 
 Let $\mcQ$ be either a $2$-nodal quartic or an irreducible quartic with one tacnode.
 Choose $N(N+1)$ contact conics $\mcC^{(i)}_j$ $i = 1, \ldots, N+1$, 
 $j = 1, \ldots, N$ to $\mcQ$ such that 
 
 \begin{itemize}
 \item $f^*_{\mcC^{(i)}_j}\mcQ$ is splitting for $1 \le j \le N - i +1$ and
 \item $f^*_{\mcC^{(i)}_j}\mcQ$ is irreducible  if otherwise.
 \end{itemize}
 
 Put 
 \[
 \mcB^{i} := \mcQ + \sum_{j=1}^N\mcC_j^{(i)} \quad i = 1, \ldots, N + 1.
 \]
 
 Then we have
 
 \begin{thm}\label{thm:z-k-plet}{
 For $i, k$ ($i \neq k$), there exists no homeomorphism
 $h : (\PP^2, \mcB^{(i)}) \to (\PP^2, \mcB^{(k)})$. In particular,
 if the combinatorial types of $\mcB^{(i)}$ $(i = 1, \ldots, N+1)$ 
 are all the same, $(\mcB^{(1)}, \ldots,\mcB^{(N)})$ is a 
 Zariski $(N+1)$-plet.
 }
 \end{thm}

 In order to show that 
 $\mcC^{(i)}_j$ ($i = 1, \ldots, N+1, j = 1, \ldots, N$) such that
 $\mcB^{(i)}$ ($i = 1, \ldots, N+1$) have the same combinatorial
 type exist, we first need to generalize Theorem~\ref{thm:q-residue}.
 
 Let $\mcQ$ be a reduced quartic and choose a smooth 
 point $z_o$ on $\mcQ$ such that 
 
 \medskip
 
 $\clubsuit$ the tangent line $l_{z_o}$ is tangent to $\mcQ$ at $z_o$ with multiplicity $2$ and 
meets $\mcQ$ at two distince residual points, or $z_o$ is an inflection point of $\mcQ$
and $I_{z_o}(l_{z_o}, \mcQ) = 3$.

\medskip

Let $\mcC$ be a contact conic to $\mcQ$. 
Then $f^*_{\mcQ}(\mcC) = \mcC^+ + \mcC^-$,  $\nu_{z_o}^*\mcC^+$  gives
rise to a section $s_{\mcC^+, z_o}$, or a bisection on 
$S_{\mcQ, z_o}$.  By Theorem~\ref{thm:shioda}, we have
a unique section for any horizontal divisor. We put
\[
s_{z_o}(\mcC^+) := \left \{ \begin{array}{cl}
                                         s_{\mcC^+, z_o} &  \mbox{if $z_o \in \mcC$}, \\
                                         \tilde{\psi}(\nu_{z_o}^*\mcC^+) & \mbox{if $z_o \not\in\mcC$.}
                                         \end{array}
                                         \right .
\]
Here $\tilde{\psi}$ is the map defined in Theorem~\ref{thm:shioda}.

 \begin{thm}\label{thm:q-residue2}{Let $\mcQ$ be a $2$-nodal quartic
 or a quartic with one tacnode.   Then
 $f_{\mcC}^*\mcQ$ is splitting if and only if $s_{z_o}(\mcC^+)$ is $2$-divisible in 
 $\MW(S_{\mcQ, z_o})$ for any $z_o$ satisfying $\clubsuit$.
 }
 \end{thm}
 
 By using Theorem~\ref{thm:q-residue2}, we can prove:
 
 \begin{thm}\label{thm:z-k-plet2}{There exist $\mcB^1, \ldots, \mcB^{N+1}$ as
 in Theorem~\ref{thm:z-k-plet} with the same combinatorial type. Moreover, there exists a Zariski $N+1$-plet of degree $2N+4$. 
 }
 \end{thm}

 In the case of $N = 2$, by combining with \cite[Theorem~1.4]{bannai16}, we have
 
 \begin{thm}\label{thm:z-5-plet}{There exists a Zariski $5$-plet of degree $8$ for
 a $2$-nodal quaritc and $2$ smooth conics.
 }
 \end{thm}
 
 \begin{rem}{\rm
 In \cite{bannai-tokunaga}, the authors constructed a Zariski 4-plet of degree 10 consisting of a quartic and three conics. 
Theorem \ref{thm:z-5-plet} gives an example of a Zariski 5-plet of degree 8. This apparently  demonstrates that the geometry of contact conics is extremely complicated but rich and worth investigating.}
 \end{rem}

 The organization of this article is as follows:
 
 In Section 1, we review our method to distinguish  the topology of $(\PP^2, \mcB)$ considered
 in \cite{bgst} and prove Theorem~\ref{thm:z-k-plet}.  The rest of this article is devoted to show
 the existence of conics $\mcC_j^{(i)}$ ($i = 1, \ldots, N+1, j = 1, \ldots, N$) such that
 $\mcB_1, \ldots, \mcB_{N+1}$ form a Zariski $(N+1)$-plet. In section 2, we introduce some
 sections of $S_{\mcQ, z_o}$ which are given by curves on $\PP^2$ not passing through $z_o$.  We prove 
 Theorem~\ref{thm:q-residue2} in Section 3. We review our method to deal with bisections
 given in \cite{bannai-tokunaga} and give examples which proves Theorem~\ref{thm:z-k-plet2} in
 Sections 4 and 5, respectively.  Theorem~\ref{thm:z-5-plet} will be proved in section 5.

\section{Proof of Theorem~\ref{thm:z-k-plet}}

%In this section, we briefly explain results, terminologies and notion which we need later. Throughout this section, let $\mcB_o$ be a (possibly empty) reduced plane curve and 
%we define  $\ucurve_{\redu}^{\mcB_o}$ to be the set of the reduced plane curves of the form
% $\mcB_o + \mcB$, where $\mcB$ is a reduced curve with no common component with $\mcB_o$.
 
 Let $\mcQ$ be an irreducible quartic and let $\mcC_1, \ldots, \mcC_N$ be contact conics to $\mcQ$.
 Put $\underline{\mcC} = \sum_{i=1}^N\mcC_i$ and for $I (\neq \emptyset) \subset \{1, \ldots, N\}$,
 put $\underline{\mcC}_I = \sum_{i \in I}\mcC_i$. Likewise \cite{bgst} we define $\Sub(\mcQ, \underline{\mcC})$ as
 \[
 \Sub(\mcQ, \underline{\mcC}):= \left \{\mcQ + \underline{\mcC}_I \mid I (\neq \emptyset) \subseteq
 \{1, \ldots, N\} \right \},
 \]
 and its subset $\Sub_k(\mcQ, \underline{\mcC})$ as
 \[
 \Sub_k(\mcQ, \underline{\mcC}):= \left \{\mcQ + \underline{\mcC}_I \mid I (\neq \emptyset) \subseteq
 \{1, \ldots, N\} , \sharp (I) = k\right \}.
 \]
 
 Now we define the map
 \begin{equation*}
	\Phi^1_{\mcQ, \mcC} : \Sub_1(\mcQ, \underline{\mcC})  \to \{0, 1\}
\end{equation*}
as follows:
\[
\Phi^1_{\mcQ, \mcC}(\mcQ + \mcC_j) =
 \left \{ \begin{array}{cl}
                                   1 & \mbox{if $f_{\mcC_j}^*\mcQ$ is splitting}\\
                                   0 &  \mbox{if $f_{\mcC_j}^*\mcQ$ is irreducible.} 
                                                            \end{array} \right .
\]

 Consider $N$ more contact conics $\mcC'_1, \ldots, \mcC'_N$ to $\mcQ$, and put
 $\underline{\mcC'} = \sum_{i=j}^N\mcC'_j$. 
 If there exists a homeomorphism $h : (\PP^2, \mcQ + \underline{\mcC}) \to 
 (\PP^2, \mcQ + \underline{\mcC'})$, then  $h$ induces a bijection
 $h_{\natural} :  \Sub(\mcQ, \underline{\mcC}) \to \Sub(\mcQ, \underline{\mcC'})$.  Moreover,
 $h_{\natural}$ induces a bijection $\Sub_1(\mcQ, \underline{\mcC}) \to \Sub_1(\mcQ, \underline{\mcC'})$
 such that
 $\Phi^1_{\mcQ, \underline{\mcC}}= \Phi^1_{\mcQ, \underline{\mcC'}}\circ h_{\natural}$: 
\[
\begin{diagram}
\node{\underline{{\rm Sub}}_1(\mcQ, \underline{\mcC})}\arrow{se,t}{\Phi^1_{\mcQ, \underline{\mcC}}}\arrow{s,l}{h_{\natural}} \\
\node{\underline{{\rm Sub}}_1(\mcQ, \underline{\mcC'})}\arrow{e,b}{ \Phi^1_{\mcQ, \underline{\mcC'}}}\node{\{0, 1\}}
\end{diagram}
\]

Under these setting, we prove Theorem~\ref{thm:z-k-plet}.
 Let  $\mcB^i= \mcQ + \underline{\mcC}^{(i)}$, $\underline{\mcC}^{(i)} = \sum_{j=1}^N\mcC_j^{(i)}$, 
 ($i = 1, \ldots, N$) be the reduced curve as in  Theorem~\ref{thm:z-k-plet}.
 Then we have 
 \[
 \sharp \left ((\Phi^1_{\mcQ, \underline{\mcC}^{(i)}})^{-1}(1)\right ) = N - i + 1.
 \]
 Hence there exists no homeomorphism $h : (\PP^2, \mcB^i) \to (\PP^2, \mcB^k)$ for $i, k$ ($i\neq k$),
 as $\mcQ$ is preserved under any homeomorphism.

 \section{Elliptic surfaces}\label{ES}

 We here summarize some facts on elliptic surfaces. We refer to
 \cite{kodaira}, \cite{miranda-basic} and  \cite{shioda90} for details.
 
 \subsection{Some general settings and facts.}
 
  Throughout this article, an elliptic surface always means a smooth projective
surface $S$ with a fibration $\varphi : S \to C$ over a smooth projective
curve, $C$, as follows:
\begin{enumerate}

\item[(i)] There exists a
{\it non empty} finite subset  $\Sing(\varphi) \subset C$ such that
$\varphi^{-1}(v)$ is a smooth curve of genus $1$ for $v \in
C\smallsetminus \Sing(\varphi)$, while $\varphi^{-1}(v)$ is not a
smooth curve of genus $1$ for $v \in \Sing(\varphi)$.

\item[(ii)] There exists a
section $O: C \to S$ (we identify $O$ with its image in $S$).

\item[(iii)]
there is no exceptional curve of the first kind in any fiber.
\end{enumerate}

For $v \in \Sing(\varphi)$, we call $F_v = \varphi^{-1}(v)$ a singular fiber over $v$. As
for  the types of singular fibers, we
use notation given by
Kodaira (\cite{kodaira}).
For $v \in \Sing(\varphi)$, we denote the irreducible  decomposition of $F_v$ by
\[
F_v = \Theta_{v, 0} + \sum_{i=1}^{m_v-1}a_{v,i}\Theta_{v,i},
\]
where $m_v$ is the number of irreducible components of $F_v$ and
$\Theta_{v,0}$ denotes the irreducible component with $\Theta_{v,0}O = 1$. We
call $\Theta_{v,0}$ {\it the identity component} of $F_v$. We also define a subset $\Red(\varphi)$
of $\Sing(\varphi)$ to be $\Red(\varphi) := \{v \in \Sing(\varphi) \mid
\mbox{$F_v$ is reducible}\}$.
For a section $s \in \MW(S)$, $s$ is said to be {\it integral} if $sO = 0$.

Let $\MW(S)$ be the set of sections of $\varphi : S \to C$. By our assumption, $\MW(S) \neq \emptyset$.
On a smooth fiber $F$ of $\varphi$, by regarding $F\cap O$ as the zero element, we can consider the abelian
group structure on $F$. Hence for $s_1, s_2 \in \MW(S)$, one can define the addition $s_1 \dot{+} s_2$ on $C\setminus
\Sing(\varphi)$. By \cite[Theorem 9.1]{kodaira}, $s_1 \dot{+} s_2$ can be extended over $C$, and we can
consider $\MW(S)$ as an abelian group. 
$\MW(S)$ is called the Mordell-Weil group.   
We also denote the multiplication-by-$m$ map ($m \in \ZZ$) on
$\MW(S)$ by $[m]s$ for $s \in \MW(S)$. Note that 
$[2]s$ is  the double of $s$ \emph{with respect to the group law on $\MW(S)$}. 
 On the other hand , we can regard the generic fiber $E:= S_{\eta}$ of
$S$ as a curve of genus $1$ over $\CC(C)$, the rational function field of $C$. The restriction of $O$
to $E$ gives rise to a ${\mathbb C}(C)$-rational point of $E$, and one can regard $E$
as an elliptic curve over ${\mathbb C}(C)$, $O$ being the zero element. By considering
the restriction to the generic fiber for each sections, $\MW(S)$ can be identified
with the set of ${\mathbb C}(C)$-rational points $E(\CC(C))$. 
For $s \in \MW(S)$, we denote the corresponding rational point by $P_s$. 
Conversely,
for an element $P \in E(\CC(C))$,  we denote the corresponding  section by $s_P$.

We also denote the addition and the multiplication-by-$m$ map on $E(\CC(C))$ by $P_1 \dot{+} P_2$ and $[m]P_1$ for
$P_1, P_2 \in E(\CC(C))$, respectively. Again, $[2]P$ is the double of $P$ \emph{with respect to the group law on $E(\CC(C))$}

Let $\NS(S)$ be the N\'eron-Severi group of $S$ and let $T_{\varphi}$ denote the subgroup of $\NS(S)$ generated by $O$, a fiber
$F$ and $\{\Theta_{v,1}, \ldots, \Theta_{v, m_v-1} \,\, (v \in \Red(\varphi))\}$. By Shioda \cite{shioda90}, we have

\begin{thm}\label{thm:shioda}(\cite[Theorem~1.2, ~1.3]{shioda90} Under our assumptions, 

\begin{enumerate}
\item[(i)] $\NS(S)$ is torsion free, and

\item[(ii)] there exists a natural map $\tilde {\psi} : \Div(S) \to \MW(S)$ which induces an isomorphism of abelian groups
\[
\psi : \NS(S)/T_{\varphi} \cong \MW(S).
\]
\end{enumerate}
\end{thm}

A curve on $S$ is said to be horizontal with respect to $\varphi$ if it does not contain any fiber components of $\varphi$. For a horizontal curve $D$ on $S$, we put $s(D):= \tilde{\psi}(D)$.
A bisection of $\varphi : S \to C$ is a reduced horizontal curve $D$ which intersects at two points with a general
fiber of $\varphi$.
As we see later, in order to construct reduced plane curve with prescribed property, we make use of a bisection $D$ and
a section $s$ with $s(D) = s$, as we see later.

% {\bf Put here summary for MWL}

In \cite{shioda90}, Shioda defined a $\QQ$-valued bilinear form $\langle \,\, , \,\, \rangle$ on $\MW(S)$ by using the intersection pairing on 
$\NS(S)$ as follows:

\begin{itemize}

\item $\langle s, \, s \rangle \ge 0$ for $\forall s \in \MW(S)$ and the equality holds if and 
only if $s$ is an element of finite order in $\MW(S)$. 

\item An explicit formula for $\langle s_1,
s_2\rangle$ ($s_1, s_2 \in \MW(S)$) is given as follows:
\[
\langle s_1, s_2 \rangle = \chi({\mathcal O}_S) + s_1O + s_2O - s_1s_2 - \sum_{v \in \Red(\varphi)}
\mbox{Contr}_v(s_1, s_2).
\]
Here we identify sections with their images, i.e., curves on $S$, and the product denotes the intersection
pairing on $\NS(S)$. Also $\mbox{Contr}_v(s_1, s_2)$ is given by
\[
\mbox{Contr}_v(s_1, s_2) = (s_1\Theta_{v,1}, \ldots, s_1\Theta_{v, m_v-1})(-A_v)^{-1}
\left ( \begin{array}{c}
        s_2\Theta_{v,1} \\
        \vdots \\
        s_2\Theta_{v, m_v-1}
        \end{array} \right ),
\]
where $A_v$ is the intersection matrix $(\Theta_{v, i}\Theta_{v, j})_{1\le i, j \le m_v - 1}$.
\end{itemize}
 As for explicit values of 
$\mbox{Contr}_v(s_1, s_2)$, we refer to \cite[(8.16)]{shioda90}.
As for explict sturcutres of $\MW(S)$ for rational ellipitic surfaces, see
\cite{oguiso-shioda}.

 \subsection{Rational elliptic surfaces, $S_{\mcQ, z_o}$}\label{sub:RES}
  
 A surface $S$ is called a rational elliptic surface if it is an elliptic surface birationally 
 equivalent to 
 $\PP^2$. 
 
 It is well-known that any rational elliptic surface is realized as a double cover
 of the Hirzebruch surface of degree $2$, which we denote by $\Sigma_2$ (see, for
example,  \cite[2.2.2]{bannai-tokunaga},
 \cite[Lecture 3, 2]{miranda-basic},  \cite[1.2]{tokunaga14}). The branch curve
 is of the form $\Delta_0 + {\mathcal T}$, where $\Delta_0$ is the section 
 with $\Delta_0^2 = -2$ and ${\mathcal T} \sim 3(\Delta_0 + {\frak f})$, ${\frak f}$
 being a fiber of the ruling of $\Sigma_2$. Note that the preimage of $\Delta_0$
 is a section $O$.
 
 On the other hand, for a reduced quartic $\mcQ$, which is not a union of four concurrent 
 and 
 a general point $z_o$ of $\mcQ$, we can associate a rational elliptic surface $S_{\mcQ, z_o}$
 as in the Introduction.

% \end{enumerate}
 
 As we have seen in \cite[Section 4]{tokunaga14}, if $z_o$ satisfies $\clubsuit$, 
 the tangent line $l_{z_o}$ at $z_o$  becomes
 an irreducible component of  a singular fiber ot type either $\I_2$ or $\III$. For other singular
 fibers, see \cite[Section 6]{miranda-persson}.
 
 Note that $\widehat{\PP^2}$ can be blown down to $\Sigma_2$ and we denote
 the composition of blowing-downs by $\bar q : \widehat{\PP^2} \to \Sigma_2$. This 
 is nothing but the realizetion of $S_{\mcQ, z_o}$ as a double cover of $\Sigma_2$.
 
 \subsection{Sections, bisections of $S_{\mcQ, z_o}$ and 
 contact conics to $\mcQ$}\label{sub:bisection}
 
 We here recall our method to treat bisections of $S_{\mcQ, z_o}$ considered in \cite[2.2.3]{bannai-tokunaga}. 
 Choose homogeneous coordinates $[T : X : Z]$ of $\PP^2$ such that $z_o = [0 : 1 :  0]$, and 
the tangent line at $z_o$ is given by $Z = 0$. $\mcQ$ is given by a homogeneous polynomial
of the form
\[
\mcQ: F(T, X, Z) = X^3Z + b_2(T, Z)X^2 + b_3(T, Z)X + b_4(T, Z) = 0,
\]
where $b_i(T, Z)$ ($i = 2, 3, 4$) are homogeneous polynomials of degree $i$. Let $U$ be an affine open set of $\PP^2$
with coordinate $(t, x) = (T/Z, X/Z)$.  Under these circumstances, the elliptic curve $E_{\mcQ, z_o}$ over $\CC(t) (\cong \CC(\PP^1))$ 
is given by the Weierstrass equation:
\[
y^2 = F(t, x, 1).
\]
Choose an element $P = (x(t), y(t)) \in E_{\mcQ, z_o}(\CC(t))$ and $r(t) \in \CC(t)$. Consider the line $L$ in ${\mathbb A}_{\CC(t)}^2$
defined by
\[
L : y = l(t, x) = r(t)(x - x(t)) + y(t).
\]
Then we have
\[
F(t, x, 1) - \{l(t, x)\}^2 = (x - x(t))g(t, x),
\]
where $g(t, x) \in \CC(t)[x]$ with $\deg_x g(t, x) = 2$. Here, $g(t,x)$ can be 
considered as a rational function on $\Sigma_2$. The zero divisor $\left(g(t,x)\right)_0$ of 
$g(t,x)$ is of the form
\[
C_{g(t,x)} + (\mbox{fiber components}),
\]
where $C_{g(t,x)}$ is a curve without any fiber components, $C_{g(t,x)}{\frak f} = 2$
and $\Delta_0 \not\subset C_{g(t,x)}$.

\begin{lem}\label{lem:bisec}(\cite[Lemma 2]{bannai-tokunaga})
{If $C_{g(t,x)}$ is irreducible and $C_{g(t,x)} \not\subset \mcT$, then

\begin{enumerate}

\item[(i)] $(q\circ f)^*C_{g(t,x)}$ is of the form
\[
(q\circ f)^*C_{g(t,x)} = C_{g(t,x)}^+ + C_{g(t,x)}^- + E,
\]
where $C_{g(t,x)}^{\pm}$ are bisections and $E$ is a divisor whose irreducible 
components are
supported on the exceptional set of $\mu$.

%all the irreducible components of $E$ are contained the exceptional set of 
%$¬?mu : S ¬?to S'$.

\item[(ii)] If we choose $C_{g(t,x)}^+$ suitably,
\[
 (y - l(t,x))|_{E_{\mcQ, z_o}} = C_{g(t,x)}^+|_{E_{\mcQ, z_o}} + P - 3 O, 
 \]
 which implies $P_{C_{g(t,x)}^+} = - P$ or equivalently  $s(C_{g(t,x)}^+) = -s_P$.

\end{enumerate}
}
\end{lem}

\begin{defin}\label{def:curve-bisec}{\rm
 For a given $P = (x(t), y(t)) \in E_{\mcQ, z_o}(\CC(t))$ and 
$r(t) \in \CC(t)$, we denote the irreducible bisection $C^+_{g(t, x)}$ obtained
in Lemma~\ref{lem:bisec} by $D(r(t), P)$.
We also denote a plane curve ${\bar q}\circ f(D(r(t), P))$ by 
$C(r(t), P)$.
}
\end{defin}

In \cite{bannai-tokunaga},
we applied Lemma~\ref{lem:bisec} to the case when 
$x(t), y(t) \in \CC[t]$ and a suitable $r(t) \in \CC[t]$ in order 
to find contact conics to $\mcQ$.

In this article, we take the same approach to find explicit example 
and to prove Theorem~\ref{thm:z-k-plet2} as that in \cite{bannai-tokunaga}.

\section{Geometricaly useful sections of $S_{\mcQ, z_o}$}

\subsection{Distinguished point free section}

Let $\mcQ$ be a reduced quartic which has at least one non-linear irreducible
component. Choose a general point $z_o$ satifying $\clubsuit$. 
Let $\varphi_{\mcQ, z_o} : S_{\mcQ, z_o} \to \PP^1$ be
the rational elliptic surface associated to $\mcQ$ and $z_o$. By our choice of $z_o$,  $S_{\mcQ, z_o}$ has a singular fiber $F_{\infty}$ 
of type $\I_2$ or $\III$, which we denote by
\[
F_{\infty} = \Theta_{\infty, 0} + \Theta_{\infty, 1},
\]
where $\Theta_{\infty, 0}$ is the exceptional divisor of the first blowing-up of $\nu : S_{\mcQ, z_o} \to S_{\mcQ}$ and 
$\Theta_{\infty, 1}$ arises from the tangent line at $z_o$, $l_{z_o}$.

Let us start with the following  lemma:

\begin{lem}\label{lem:line-sec}{ Let $s \in \MW(S_{\mcQ, z_o})$ such that $sO = 0$ and $s\Theta_{\infty, 1} = 1$. Then
${\bar q}\circ f(s)$ is a line $L_s$ in $\PP^2$ such that
\begin{enumerate}

\item[(i)] $z_o \not\in L_s$, and

\item[(ii)] $I_x(L_s, \mcQ)$ is even for $\forall x \in L_s\cap \mcQ$.

\end{enumerate}

Conversely, any line satisfying the above two conditions gives rise to sections $s_{L^{\pm}}$ such that
$s_{L^{\pm}} O = 0$ and $s_{L^{\pm}} \Theta_{v, 1} = 1$.
}
\end{lem}

\proof By our construction, ${\bar q}\circ f(s)$ intersects any line through $z_o$ at one point. Hence ${\bar q}\circ f(s)$ is
a line not passing through $z_o$. Also, if $I_x(L_s, \mcQ)$ is odd for some $x \in L_s\cap \mcQ$,
$({\bar q}\circ f)^*(L_s)$ is a sum of an irreducible curve which is mapped to $L_s$ and exceptional curves for $\mu$.
On the other hand, both $s$ and $[-1]_{\varphi_{\mcQ, z_o}}s$ are contained in $({\bar q}\circ f)^*L_s$, which leads
us to a contradiction. 

Conversely, if $L$ is a line satisfying the conditions (i) and (ii), then $(f'\circ \mu)^*L$ is of 
the form
\[
({\bar q}\circ f)^*L = L^+ + L^- + E,
\]
where $E$ consist of exceptional curves for $\mu : S_{\mcQ} \to S'_{\mcQ}$. Since $L^{\pm}$
meet $(f'\circ \mu)^*l$ ($l$ is a line through $z_o$) at one point. Hence
$\nu^*L^{\pm}$ give rise to  sections, which we denote by $s_{L^{\pm}}$. As $z_o \not\in L$,
$s_{L^{\pm}} O = 0$ and $s_{L^{\pm}}\Theta_{\infty, 1} = 1$.
\proofend

\begin{defin}\label{def:dp-free-sec}{\rm
For $s \in \MW(S_{\mcQ, z_o})$, we call $s$ a {\it distinguished point free section (dp-free section, for
short)}  if
$sO = 0$ and $s\Theta_{\infty, 1} = 1$. Otherwise, we call $s$ a {\it distinguished point
sensitive section (dp-sensitve section, for short)}.
}
\end{defin}

\begin{rem}\label{rem:dp-free-sec}{\rm
By Lemma~\ref{lem:line-sec}, any distinguished point free section of $S_{\mcQ, z_o}$ arises
from a line in $\PP^2$ satisfying the conditions (i) and (ii).
}
\end{rem}

 \subsection{Examples}
 
 \begin{exmple}\label{eg:smooth}{\rm (\cite{shioda93})
  For a smooth quartic curve $\mcQ$,
 it is well-known that $\mcQ$ has $28$ bitangent lines, which give rise to
 $56$ dp-free sections, which generate $\MW(S_{\mcQ, z_o})$.  In \cite{shioda93}, these $56$ sections are intensively
 studied  and their applications are given. See \cite{shioda93}, for detail.
 }
 \end{exmple}
 
 \begin{exmple}\label{eg:2-nodal}{\rm
 Let $\mcQ$ be an irreducible quartic  with $2$ nodes only as its singularities. Choose a
 smooth point $z_o$ satisfying the condition $(\clubsuit)$. In this case, The configuration of reducible
  singular fibers of $\varphi_{\mcQ, z_o} : S_{\mcQ, z_o} \to \PP^1$ is either
 $3\I_2$ or $2\I_2, \III$. For both cases, $\MW(S_{\mcQ, z_o}) \cong A_1^*\oplus D_4^*$.
 We consider a basis of $\MW(S_{\mcQ, z_o})$ consisting of dp-free sections.
 
 Let $L_0$ be the line connecting $x_1$ and $x_2$.  $L_0$ gives rise to $2$ dp-free sections $s^+_{0}$ and $s^-_{0}$.
 
 By considering the projection centered at $x_i$ ($i = 1, 2$), there exist $4$ tangent lines
 $L_{i, j}$ ($j = 1, 2, 3, 4$) through $x_i$. Here, the case when $L_{i,j}$ is tangent to one of the branches of 
 the node with multiplicity $3$ is also considered.  
 Since
 \[
 (q\circ f_{\mcQ})^*L_{i, j} = L_{i, j}^+ + L_{i,j}^-,
 \]
 these $8$ lines give rise to $16$ dp-free sections $s^+_{i, j} = \nu_{z_o}^*L^+_{i, j}, 
 s^-_{i, j} = \nu_{z_o}^*L^-_{i, j}$
 ($i = 1, 2$, $j = 1, 2, 3, 4$). By (re)labeling $i, j, \pm$ suitably, we may assume that
 \[
 s^+_{1, k}s^+_{1, l} = -\delta_{kl}, \quad s^+_{1, j}s^+_{2, 1} = 1 \quad (j = 1, 2, 3, 4),  \quad s^\pm_{0}s^{\pm}_{i, j} = 0 \quad (i = 1,2,\, j =1, 2, 3, 4)
 \]
 Hence, we have
 \[
 \langle s^+_{1,k}, s^+_{1,l} \rangle = \delta_{kl}, \,\,   \langle s^+_{1, k}, s^+_{2,1}\rangle =  - \frac 12, \,\,  \langle s^+_0, s^+_0 \rangle = \frac 12, \,\,  \langle s^+_0, s^+_{i, j} \rangle = 0 \,\,  (i =1, 2, j = 1,2, 3, 4)
 \]
 Put
 \[
 s_{z_o, 0} = s^+_0 = \nu_{z_o}^*L_0^+, \,\, s_{z_0, i} = s^+_{1,i} = \nu_{z_o}^*L^+_{1,i},\, (i =1, 2, 3) \,\,  
  s_{z_o,4} = s^+_{2,1} = \nu_{z_o}^*L_{2,1}^+.
 \]
 Then we have the Gram matrix $\left [ \langle s_{z_o, i}, s_{z_o, j} \rangle \right ]$
 \[
 \left [ \begin{array}{ccccc}
        \frac 12 & 0 & 0 & 0 & 0 \\
        0 & 1 & 0 & 0 & -\frac 12 \\
        0 & 0 & 1 & 0 & -\frac 12 \\
        0 & 0 & 0& 1 & -\frac 12 \\
        0  & -\frac 12 & - \frac 12 & - \frac 12 & 1
        \end{array}
        \right ]. 
 \]
 Let $M$ be a sublattice of $\MW(S_{\mcQ, z_o})$ given by $s_{z_o, 0}, 
 s_{z_o, 1}, s_{z_o, 2}, s_{z_o, 3}, s_{z_o, 4}$.  As $\det M = 1/8$. $\MW(S_{\mcQ, z_o})
 = M$. Finally if we put
$ t_0  =  s_{z_o, 0}$,
 $t_1  =  s_{z_o, 1}$,
 $t_2  =  s_{z_o, 1} + s_{z_o, 2}$, 
 $t_3  =  -(s_{z_o, 3} + s_{z_o, 4})$,
 $t_4  =  - s_{z_o, 4}$,
% Then we have
% \[
% \left [\langle t_i, t_j \rangle \right ] = \left [ \begin{array}{ccccc}
%        \frac 12 & 0 & 0 & 0 & 0 \\
%        0 & 1 & 1 & \frac 12 & \frac 12 \\
%        0 & 1 & 2 & 1 & 1 \\
%        0 & \frac 12 & 1& 1 & \frac 12 \\
%        0  & \frac 12  & 1 & \frac 12 &  1
%        \end{array}
%        \right ]. 
% \]
 we get a basis for $A_1^*\oplus D_4^*$ with the well known Gram matrix.
  }
 \end{exmple}
 
 \begin{exmple}\label{eg:tacnode}{\rm
 Let $\mcQ$ be an irreducible quartic with one tacnode, $x_o$. Choose $z_o \in \mcQ$ satisfying $\clubsuit$ and let $\varphi_{\mcQ, z_o} : S_{\mcQ, z_o} \to \PP^1$ be
 the rational elliptic surface as in the Introduction.   The configuration of reducible singular fibers of $\varphi_{\mcQ, z_o}$ is either  $\I_4, \, \I_2$ or $\I_4, \,  \III$. By \cite{oguiso-shioda},  
 $\MW(S_{\mcQ, z_o}) \cong A_1^* \oplus A_3^*$. We also consider a basis consisting of
 do-free sections as in Example~\ref{eg:2-nodal}.
 
  Let $L_{x_o}$ be the tangent
 line at $x_o$.  By  a similar argument to the one in Example~\ref{eg:2-nodal}, there exist
 $4$ lines  $L_i$ ($i = 1, \ldots, 4$) through $x_o$ and tangent to $\mcQ$ at another residual 
 point. By our construction of $S_{\mcQ}$, we have
 \[
 (q\circ f_{\mcQ})^*(L_{x_o}) = L^+_{x_o} + L^-_{x_o} + E_{x_o},  \quad 
 (q\circ f_{\mcQ})^*(L_{i}) = L^+_{i} + L^-_{i} + E_{i}, \, \, (i = 1, \ldots, 4),
 \]
 where $E_{x_o}$, $E_i$ ($i = 1, \ldots, 4$) are divisors whose supports 
 are contained in the exceptional set of $\mu$.  By Lemma~\ref{lem:line-sec},  we infer that 
 $\nu_{z_o}^*L_{x_o}^{\pm}$, $\nu_{z_o}^*L_i^{\pm}$ $(i = 1, \ldots, 4)$ are 
 dp-free sections of $\varphi_{\mcQ, z_o}$.
 Put 
 \[
 s_{z_o, 0}:= \nu_{z_o}^*L_{x_o}^+, \quad s_{z_o, i}:= \nu_{z_o}^*L_i^+.
 \]
 By labelling appropriately, we may assume:
 \[
 \begin{array}{cl}
 \langle s_{z_o, 0}, s_{z_o, 0}\rangle = \frac 12, & 
 \langle s_{z_o,0}, s_{z_o, i} \rangle = 0, \,\, (i = 1, \ldots, 4) \\
 \langle s_{z_o, i}, s_{z_o, i}\rangle = \frac 34, & \langle s_{z_o, i}, s_{z_o, j}\rangle = -\frac 14, \,
 (i\neq j)\,\,
  (i, j = 1, \ldots, 4)
  \end{array}
  \]
 Hence the Gram matrix for $s_{z_o, i}$ $(i = 0, \ldots 3)$ is 
 \[
 \left [
 \begin{array}{cccc}
 \frac 12 & 0 & 0 & 0 \\
 0 & \frac 34 & -\frac{1}4 & -\frac {1}4 \\
 0 & -\frac{1}4 & \frac 34 & -\frac {1}4 \\
 0 & -\frac{1}4 & -\frac{1}4 & \frac 34 
 \end{array}
 \right ]
 \]
 Since the determinant of the above matrix is $1/8$, we infer that $s_{z_o, i}$ $(i = 0, \ldots, 3)$ form
 a basis of $\MW(S_{\mcQ,z_o})$. Now put $t_0 := s_{z_o, 0}, \, t_1:= s_{z_o, 1}, \,  t_2:= s_{z_o, 1} + 
 s_{z_o, 2}, \, t_3 := -s_{z_o,3}$. The Gram matrix $[\langle t_i,t_j\rangle]$ becomes the well known one.
% \[
% \left [
% \begin{array}{cccc}
% \frac 12 & 0 & 0 & 0 \\
% 0 & \frac 34 & \frac 12 & \frac 14 \\
% 0 & \frac 12 & 1 & \frac 14 \\
% 0& \frac 14 & \frac 12 & \frac 34
% \end{array}
% \right ].
% \]
% basis of $A_1^*\oplus A_3^*$.
 }
 \end{exmple}
 
 \begin{exmple}\label{eg:3-nodal}{\rm
 Let $\mcQ$ be an irreducible quartic with $3$ nodes, which we denote $x_1, x_2, x_3$. Choose $z_o$ satisfying
 $\clubsuit$ in the Introduction. In this case, the configuration of singular fibers of $\varphi_{\mcQ, z_o} : S_{\mcQ, z_o}
 \to \PP^1$ is either $4\I_2$ or $3\I_2, \III$. For both cases, $\MW(S_{\mcQ, z_o}) \cong (A_1^*)^{\oplus 4}$. As we have
 seen in \cite{tumenbayar-tokunaga}, a line connecting $2$ nodes $x_i$ and $x_j$ gives rise to two dp-free sections $s_{i, j}^{\pm}$,
 which are generators
 of one of the direct summands of $\MW(S_{\mcQ, z_o})$.  
 Hence $3$ direct summands are generated these
 dp-free sections which arise from lines connecting $2$ nodes. We then choose a line $l_4$ which passes through $x_3$ and
 tangent to $\mcQ$ another point or  passes through $x_3$ with multiplicity $4$. Then $l_4$ gives another dp-free section
 $s_4^{\pm}$. By relabelling $s_{i, j}^{\pm}$ suitably, we may assume that
 \[
 s_4^+ s_{i, j}^+ = 0.
 \]
 Put $s_1:= s_{1, 2}^+, s_2:= s_{1,3}^+, s_3:= s_{2, 3}^+, s_4:= s_4^+$, and we have the Gram matrix
 \[
 [\langle s_i, s_j \rangle] 
 = \left [ \begin{array}{cccc}
            \frac 12 & 0 & 0 & \frac 12 \\
              0 & \frac 12 & 0 & 0 \\
              0 &  0 & \frac 12 & 0 \\
              \frac 12 & 0 & 0 & 1
            \end{array}
            \right ]  
 \]
 Since the determinant of the above matrix is $1/16$, we infer that $s_1, \ldots, s_4$  form a basis of $\MW(S_{\mcQ, z_o})$.

\begin{rem}
It is an interesting problem to determine which rational elliptic surfaces admit a dp-free basis, and also to find an explicit geometric description of each basis.
\end{rem}
% The remaining direct summand is generated
% by sections $s_0^{\pm}$arising from a unique conic tangent at $z_o$ passing through $3$ nodes.
% Thus we have basis $\{s_0^+, s_{1,2}^+, s_{2, 3}^+, s_{1,3}^{\pm}\}$ of $\MW(S_{\mcQ, z_o})$. In this case, 
% it is impossible to find a basis consisting of dp-free sections. 
 
 }
 \end{exmple}

\section{Proof of Theorem~\ref{thm:q-residue2}}\label{sec:proofthm:q2}

% Let $\mcQ$ be a $2$-nodal quartic and let $\mcC$ be a contact conic to $\mcQ$. Let
% $f_{\mcC} : Z_{\mcC} \to \PP^2$ be the double cover of $\PP^2$ branched along
% $\mcC$ and we denote the covering transformation of $f_{\mcC}$ by
% $\sigma_{f_{\mcC}}$. The pullback $f_{\mcC}^*\mcQ$ is eihter
% irreducible, or splitting, i.e., of the form $\mcQ^+ + \mcQ^-$, $\sigma_{f_{\mcC}}^*\mcQ^+ = 
% \mcQ^-$. As we have seen in \cite{tokunaga10}, a criterion for $f_{\mcC}^*\mcQ$ to be
% irreducible or not can be formulated as follows:
% 
% Choose a point $z_o \in \mcQ\cap \mcC$ and consider the rational elliptic surface
% $\varphi_{\mcQ, z_o} : S_{\mcQ, z_o} \to \PP^1$ as in previous sections. 
% $(q\circ f_{\mcQ})^*\mcC = \mcC^+ + \mcC^-$ and $\nu_{z_o}^*\mcC^{\pm}$ 
% are sections. We denote them by $s_{\mcC^{\pm},z_o}$, respectively. Then we have
 
% \begin{thm}\label{thm:q-residue}{\cite[Theorem~2.1]{tokunaga10} $f_{\mcC}^*\mcQ$ is
% splitting if and only if $s_{\mcC^+, z_o}$ is $2$-divisible in $\MW(S_{\mcQ, z_o})$
% }
% \end{thm}
% 
% Note that $s_{\mcC^-, z_o}$ is $2$-divisible if and only if so is $s_{\mcC^+, z_o}$,
% as $s_{\mcC^-, z_o} = [-1]s_{\mcC^+, z_o}$.
 
\subsection{An application of distinguished point free sections}

  We first consider an application of distinguished point free sections. Let $\mcQ$ be as before.
 Let $f'_{\mcQ} : S'_{\mcQ} \to \PP^2$ be the double cover branched along $\mcQ$ and 
 $S_{\mcQ}$ the canonical resolution of the singularities of $S'_{\mcQ}$. Choose two distinguished 
 point $z_1$ and $z_2$ satisfying $\clubsuit$.
 The resolution maps for the pencils $\Lambda_{z_1}$ and $\Lambda_{z_2}$
 will be denoted by $\nu_i:S_{\mcQ, z_i}\rightarrow S_{\mcQ}$, respectively.  Each $\nu_i$ is a composition of
 two blowing-ups. 

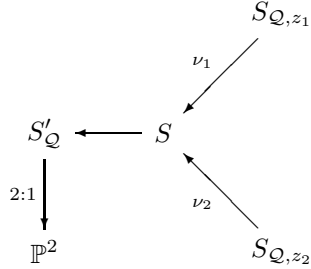
\[
\begin{diagram}
\node{}\node{}\node{S_{\mcQ, z_1}}\arrow{sw,t}{\nu_1}\\
\node{S'_{\mcQ}}\arrow[1]{s,l}{2:1}\node{S}\arrow{w}\node{}\\
\node{\PP^2}\node{}\node{S_{\mcQ, z_2}}\arrow{nw,b}{\nu_2}
\end{diagram}
\]

The exceptional divisor of the second blowing-up of $\nu_1$ (resp. $\nu_2$) gives rise to a section of 
$S_{\mcQ, z_1}$ (resp. $S_{\mcQ, z_2}$). We will regard this section as the zero section and denote 
it by $O_{1}$ (resp. $O_2$). By construction, $S_{\mcQ, z_1}$ and $S_{\mcQ, z_2}$ are rational elliptic 
surfaces that have the  same configuration of singular fibers, except possibly the one arising from
tangent lines $l_{z_i}$ at $z_i$ which is either   of type $\I_2$ or $\III$.  

Let $D_1,\ldots,D_m$ be divisors on $S$ such that they do not pass through  
$(f'_{\mcQ}\circ \nu_1)^{-1}(z_1), (f'_{\mcQ}\circ\nu_2)^{-1}(z_2)$ and 
their strict transforms under $\nu_1$ (resp.  $\nu_2$) give rise to sections of 
$S_{\mcQ, z_1}$ (resp. $S_{\mcQ, z_2}$). Note that these sections are distinguished point free sections. Let $s_i(D_j)$ denote the 
section corresponding to $D_j$ on $S_{\mcQ, z_i}$. 
Put $C_i=s_i(D_1)\dot+\cdots\dot+s_i(D_m) \in \MW(S_{\mcQ, z_i})$. 
Put $\overline{C_i} = \nu_i(C_i)$ and let $\widehat{C}_i$ be the strict 
transform of $\overline{C}_i$ under $\nu_j^{-1}$.  $(i\not=j)$. For $z_1$ and $z_2$ with $\clubsuit$,
 $\widehat{C}_2$ becomes a multi-section of $S_{z_1}$. 
 Under this setting, we have:

\begin{thm}\label{thm:prescription}
{
Suppose that $C_2\not=O_2$ and $z_1\not\in \overline{C_2}$. Then
\[
s(\widehat{C}_2)=C_1.
\]
%where $\bar{\psi}_1$ is the map from $\NS(S_{p_1})$ to $\MW(S)$ given in \cite{shioda90}.
}
\end{thm}

\proof
Note that each surface $S_{\mcQ, z_i}$ 
has a $\I_2$ or $\III$ type singular fiber whose components arise from the exceptional divisor of the first 
blow up in $\nu_i$ which meets $O_{i}$ and the strict transform of the tangent line $l_{z_i}$ of $\mcQ$ at $z_i$. We will denote these components by $\Theta_{z_i,0}$ and $\Theta_{z_i,1}$.  All the other reducible singular fibers arise from the exceptional sets of the resolution $S_{\mcQ} \rightarrow S_{\mcQ}^\prime$, 
hence they are in 1 to 1 correspondence with the singularities of $Q$. We will denote their components by $\Theta_{v,i}$ where $v\in \Sing(\mcQ)$.
Let $D^\prime=s_1(D_1)+\cdots+s_1(D_m)$. Note that  the sum taken here is regarded as a sum of divisors on $S_{\mcQ, z_1}$. Then since the Abel-Jacobi map $\tilde{\psi}_i$ for each surface is a homomorphism, $\tilde{\psi}_1(D^\prime)=s_1(D_1)\dot+\cdots\dot+s_1(D_n)=C_1$. 
Hence by \cite[Lemma~5.1]{shioda90} we have the equivalence
\[
D^\prime\underset{S_{\mcQ, z_1}}\sim (C_1)+d(O_1)+nF_1+a\Theta_{z_1,1}+\sum_{v\in \Sing(\mcQ)}\sum_{i=1}^{m_v-1} b_{v,i}\Theta_{v,i},
\]
where $\underset{S_{\mcQ, z_1}}\sim$ denotes linear equivalence of divisors on $S_{\mcQ, z_1}$
 
 Similarly for $D^{\prime\prime}=s_2(D_1)+\cdots+s_2(D_m)$, we have the equivalence
\[
D^{\prime\prime}\underset{S_{\mcQ, z_2}}\sim (C_2)+d(O_2)+nF_2+a\Theta_{z_2,1}+\sum_{v\in \Sing(\mcQ)}\sum_{i=1}^{m_v-1} b_{v,i}\Theta_{v,i}.
\]
Note that by construction of $D^\prime, D^{\prime\prime}$ and \cite[Theorem~9.1]{kodaira} (or 
\cite[1.1]{tokunaga12}, 
the coefficients $d, n, a, b_{v,i}$ are the same in both cases. Since ${\nu_1}_\ast(D^\prime)=D_1+\cdots+D_m={\nu_2}_\ast(D^{\prime\prime})$, by \cite[Thorem~1.4]{fulton}, we have 
\begin{align*}
\bar{C}_1+&nF_1+a\nu_{\ast}\Theta_{z_1,1}+\sum_{v\in \Sing(Q)}\sum_{i=1}^{m_v-1} b_{v,i}\Theta_{v,i}\\
&\underset{S_{\mcQ}}{\sim} \bar{C_2}+nF_2+a\nu_{\ast}\Theta_{z_2,1}+\sum_{v\in \Sing(Q)}\sum_{i=1}^{m_v-1} b_{v,i}\Theta_{v,i}.
\end{align*}

 Then since $F_1\underset{S_{\mcQ}}{\sim}F_2$ and $\Theta_{z_1,1}\underset{S_{\mcQ}}{\sim}\Theta_{z_2,1}$, because they are inverse images of lines of $\PP^2$, we obtain the equivalence
\[
\bar{C}_1\underset{S_{\mcQ}}{\sim}\bar{C}_2.
\]
By pulling this equivalence back by ${\nu_1}$, we obtain
\begin{align*}
\hat{C}_2&\underset{S_{\mcQ, z_1}}\sim C_1+\alpha O_1+\beta\Theta_{z_1,0}\\
&\underset{S_{\mcQ, z_1}}\sim C_1+\alpha O_1+\beta F-\beta \Theta_{z_1,1},
\end{align*}
for some integers $\alpha$ and $\beta$. Hence by Theorem~\ref{thm:shioda} we have $\tilde{\psi}_1(\hat{C}_2)=C_1$.
\qed

\subsection{Proof of  Theorem~\ref{thm:q-residue2}}

Let $\mcQ$ be a quartic and $z_o\in\mcQ$ be a point satisfying $\clubsuit$.  Furthermore, 
assume that $\MW(S_{\mcQ, z_o})$ is generated by dp-free sections $s_{z_o,i}$ $(i=0,\ldots, k)$ as in the setting of the previous subsection.
% Let $s_{z_o, i}$ $(i = 0, \ldots, 4)$ be elements of $\MW(S_{\mcQ, z_o})$ as in
% Example~\ref{eg:2-nodal}. Note that these sections form a basis of
% $\MW(S_{\mcQ, z_o})$ for $z_o$ satisfying $\clubsuit$ as we see in 
% Example~\ref{eg:2-nodal}. 
Let $\mcC'$ be a contact conic with $z_o \not\in \mcC'$.
 On $S_{\mcQ}$,  $(q\circ f_{\mcQ})^*\mcC^\prime = {\mcC'}^+ + {\mcC'}^-$, and
 $\nu^*_{z_o}{\mcC'}^{\pm}$ become bisections on $S_{\mcQ, z_o}$. 
 Let $s_{z_o}(\nu^*_{z_o}{\mcC'}^+) (= \tilde{\psi}((\nu^*_{z_o}{\mcC'}^+))$ be the section 
 determined by ${\mcC'}^+$. Assume that
 \[
 s_{z_o}(\nu^*_{z_o}{\mcC'}^+)  = \sum_{i=0}^k [a_i]s_{z_o, i}, \quad a_i \in \ZZ.
 \]
 Choose $z'_o \in \mcC'\cap \mcQ$. Then $\nu_{z'_o}^*{\mcC'}^+$ is an element of
 $\MW(S_{\mcQ, z'_o})$, which we denote by $s_{z'_o, {\mcC'}^+}$. Then we have
 
 \begin{prop}\label{prop:key-split}{Under the above setting,
 \[
 s_{z'_o, {\mcC'}^+} = \sum_{i=0}^k[a_i]s_{z'_o, i}.
 \]
 }
 \end{prop}
 
 \proof Put $s_{z'_o, {\mcC'}^+} = \sum_{i=0}^4[b_i]s_{z'_o, i}$. By Theorem~\ref{thm:prescription},
 we have
  \[
 s_{z_o}(\nu^*_{z_o}{\mcC'}^+)  = \sum_{i=0}^4 [b_i]s_{z_o, i}, \quad b_i \in \ZZ.
 \]
By our assumption, we have $a_i = b_i$ $(i = 0, \ldots, k)$.
\proofend

{\bf Case I: $\mcQ$ is a $2$-nodal quartic.}

By Example \ref{eg:2-nodal}, $\MW(\mcQ_{z_o})$ is generated by bp-free sections. Furthermore, since $\langle s_{z'_o, {\mcC'}^+}, s_{z'_o, {\mcC'}^+} \rangle = 2$ and 
$\MW(S_{\mcQ, z'_o}) \cong A_1^*\oplus D_4^*$,  we have

\begin{cor}\label{cor:key-split}{ $f_{\mcC'}^*\mcQ$ is split if and only if
$a_0  = 2, a_1= \ldots a_4 = 0$.
}
\end{cor}

Also by \cite[Theorem~0.2]{tokunaga10},  we have

\begin{cor}\label{cor:key-split2}{There exists a Galois cover of $\PP^2$ such that the Galois group
is isomorphic to the dihedral group of oder $2p$ and the branch locus is $\mcC + \mcQ$ if 
and only if $a_0  = 2, a_1= \cdots=a_4 = 0$.
}
\end{cor}

%\subsection{Proof of Theorem~\ref{thm:q-residue2} for  $\mcQ$: an irreducible quartic with one tacnode}
%
% 
% 
% Let $s_{z_o, i}$ $(i = 0, \ldots, 3)$ be elements of $\MW(S_{\mcQ, z_o})$ as in
% Example~\ref{eg:tacnode}. Note that these sections form a basis of
% $\MW(S_{\mcQ, z_o})$ for $z_o$ satisfying $\clubsuit$ as we see in 
% Example~\ref{eg:tacnode}. 
% Let $\mcC'$ be a contact conic with $z_o \not\in \mcC'$. Define
% $s_{z_o}(\nu^*_{z_o}{\mcC'}^+)$ in the same way as in the case of $\mcQ$: $2$-nodal and 
% assume that
% \[
% s_{z_o}(\nu^*_{z_o}{\mcC'}^+)  = \sum_{i=0}^ 3[a_i]s_{z_o, i}, \quad a_i \in \ZZ.
% \]
% Choose $z'_o \in \mcC'\cap \mcQ$. $\nu_{z'_o}^*{\mcC'}^+$ is an element of
% $\MW(S_{\mcQ, z'_o})$, which we denote by $s_{z'_o, {\mcC'}^+}$.  The similar
% argument to that of $2$-nodal case implies the following:
% 
% 
%  \begin{prop}\label{prop:key-split2}{
% \[
% s_{z'_o, {\mcC'}^+} = \sum_{i=0}^3[a_i]s_{z'_o, i}.
% \]
% }
% \end{prop}
% 
% Also we have the following corollaries:
 {\bf Case II: $\mcQ$ is an irreducible quartic with a tacnode.}

By Example \ref{eg:tacnode}, $\MW(\mcQ_{z_o})$ is generated by bp-free sections. Furthermore, since $\langle s_{z'_o, {\mcC'}^+}, s_{z'_o, {\mcC'}^+} \rangle = 2$ and 
$\MW(S_{\mcQ, z'_o}) \cong A_1^*\oplus A_3^*$,  we have
 \begin{cor}\label{cor:key-split3}{ $f_{\mcC'}^*\mcQ$ is split if and only if
$a_0  = 2, a_1= a_2 =  a_3 = 0$.
}
\end{cor}

Also by \cite[Theorem~0.2]{tokunaga10},  we have

\begin{cor}\label{cor:key-split4}{There exists a Galois cover of $\PP^2$ such that the Galois group
is isomorphic to the dihedral group of oder $2p$ and the branch locus is $\mcC + \mcQ$ if 
and only if $a_0  = 2, a_1= a_2 =  a_3 = 0$.
}
\end{cor}

The above Corollaries combined proves Theorem \ref{thm:q-residue2}.

\section{Proof of Theorem~\ref{thm:z-k-plet2}: Existence of contact conics}\label{sec:proof_thm:z-k-plet2} 

In this section we prove Theorem~\ref{thm:z-k-plet2} by explicitly constructing the desired contact conics. First, we describe the method of constructing contact conics in general, and afterwords give explicit equations.

%
%\subsection{General method of constructing  contact conics}\label{sub:conic-construction}
%
We utilise the method to construct bisections described in  Section~\ref{sub:bisection}.
We assume that $S_{\mcQ, z_o}$ is given the Weierstrass equation as in  section~\ref{sub:bisection}.

\bigskip

{\bf Case I: $\mcQ$ is a $2$-nodal quartic.}

Let $P_i$ $(i =0, \ldots, 4)$ be points in $E_{\mcQ, z_o}(\CC(t))$ corresponding to
the sections $s_{z_o, i}$ $(i = 0, \ldots, 4)$ in Example~\ref{eg:2-nodal}. 
Put  $[2]P_0 = (x_1(t), y_1(t))$, $P_1\dot{+}P_2 = (x_2(t), y_2(t))$.  We apply the method to construct
bisections described in section~\ref{sub:bisection} as follows:

\begin{itemize} 

\item We choose $r_{i, a}(t) \in \CC(t), a \in \CC$ $(i = 1, 2)$ appropriately such that both
$C(r_{a, 1}, [2]P_0)$ and $C(r_{a, 2}, P_1\dot{+}P_2)$ give conics.

\item Choose $a_1, \ldots, a_N, b_1, \ldots, b_N \in \CC$ such that  
$C(r_{a_i, 1}, [2]P_0)$  and $C(r_{b_i, 2}, P_1\dot{+}P_2)$ are contact conics to $\mcQ$
contact at $4$ distinct points.

\item The $2N$ conics as above meets transversely with each other.

\end{itemize}

Put  $\mcC_i = C(r_{a_i, 1}, [2]P_0)$ $(i = 1, \ldots, N)$ and $\mcC'_j := C(r_{b_j, 2} P_1\dot{+} P_2)$
$(j = 1, \dots, N)$. Then we have
\begin{lem}\label{lem:split}{Under the  setting given above, 
\begin{enumerate}
\item[(i)] $f_{\mcC_i}^*\mcQ$ is splitting.

\item[(ii)] $f_{\mcC'_j}^*\mcQ$ is irreducible.
\end{enumerate}
}
\end{lem}

\proof  (i) Choose $z'_o \in \mcQ\cap \mcC_i$. 
Note that $\mcC_i$ is of the form $\mcC^+ + \mcC^-$ on $S_{\mcQ}$. By Proposition~\ref{prop:key-split}, 
we may assume that  
$\nu_{z'_o}^*C^+ = [2]s_{z'_o, 0}$. Hence by Theorem~\ref{thm:q-residue}, $f_{\mcC_i}^*\mcQ$ is 
splitting.

(ii)  Choose $z'_o \in \mcQ\cap \mcC'_j$. 
Again, $\mcC_i$ is of the form $\mcC^+ + \mcC^-$ on $S_{\mcQ}$. By Proposition~\ref{prop:key-split}, we
may assume that
$\nu_{z'_o}^*C^+ = s_{z'_o, 1} \dot+ s_{z'_o,2}$. Hence by Theorem~\ref{thm:q-residue}, $f_{\mcC'_j}^*\mcQ$ is 
irreducible.
\proofend

Now in order to show the existence of conics as in Theorem~\ref{thm:z-k-plet2}, 
we choose $N$-conics from the above $2N$ conics suitably.

\bigskip

{\bf Case II:  $\mcQ$ is an irreducible quartic with one tacnode.}

Let $P_i$ $(i =0, \ldots, 3)$ be points in $E_{\mcQ, z_o}(\CC(t))$ corresponding to
the sections $s_{z_o, i}$ $(i = 0, \ldots, 4)$ in Example~\ref{eg:tacnode}. 
Put  $[2]P_0 = (x_1(t), y_1(t))$, $P_1\dot{-}P_2 = (x_2(t), y_2(t))$.  We apply the method to construct
bisections described in section~\ref{sub:bisection}. The remaining argument is almost the same
as Case I, and we omit it.

%We construct contact conics to $\mcQ$ via the method described in
%section~\label{sub:conic-construction}, by which we have Theorem~\ref{thm:z-k-plet2}.
%

\bigskip

\subsection{Case I: $\mcQ$ is a $2$-nodal quartic}

Let $F(T, X, Z)$ be a homogeneous polynomial 
 \[
 {X}^{3}Z+ \left( 271350Z-98\,T \right) {X}^{2}+T\left( T-5825Z \right) 
 \left( T-2025Z \right)X+36\,{T}^{2} \left(T-2025Z \right)^{2}
 \]
and let $\mcQ$ be a quartic given by $F = 0$. $\mcQ$ is a $2$-nodal quartic and it has
two nodes at $x_1=[0:0:1]$ and $x_2=[2025:0:1]$. Let $z_o=[0:1:0]$. 
%Let $S$ be the rational elliptic surface whose generic fiber is given by the Weierstrass equation
%\begin{align*}
% y^2&= F(t, x, 1) \\
%  & ={x}^{3}+ \left( 271350-98\,t \right) {x}^{2}+t \left( t-5825 \right) 
% \left( t-2025 \right) x+36\,{t}^{2} \left( t-2025 \right) ^{2}.
% \end{align*}
 The associated rational elliptic surface $S_{\mcQ,z_o}$ was given and studied by Shioda and Usui in \cite[p. 198]{shioda-usui}.
% it corresponds to considering the rational elliptic surface associated to the quartic $\mathcal{Q}$
%  and $z_o=[0:1:0]$. 
 According to \cite{shioda-usui},
  $S_{\mcQ,z_o}$ has two singular fibers of type $\I_2$ and one singular fiber of type $\III$. The Mordell-Weil lattice of $S$ is  $\MW(S)\cong A_1^\ast\oplus D_4^\ast$ and the narrow Mordell-Weil lattice is  $\MW(S)^0\cong A_1\oplus D_4$. 
 
  Let $L_0, L_{1,1}, L_{1,2}, L_{1,3}, L_{2,1}$ be lines defined by
 \begin{align*}
 L_0:\, X=0, 
\quad L_{1,1}:  \,32T+X=0, 
\quad L_{1,2}:  \,28T-X=0,\\ 
 L_{1,3}:  \,20T+X=0, 
 \quad L_{2,1}:  \,35T+X-70875Z=0.
 \end{align*}
 The lines $L_{i,j}$ pass through $x_i$ and is tangent to $\mcQ$.
 
 The lines give rise to dp-free sections $s_0, \ldots, s_4$ with coordinates
 \begin{align*}
 s_0=(0,6t(t-2025)), 
 \quad s_1=(-32t, 2t(t-3465)), 
\quad s_2=(28t, 8t(t-3285)), \\
 s_3=(-20t, 4t(t-1125)),
\quad s_4=(-35t+70875, (t+20475)(t-2025)),
 \end{align*}
% Let $s_{z_o,0}=s_0^+, s_{z_o,1}=s_{1,1}^+, s_{z_o,2}=s_{1,2}^+, s_{z_o, 3}=s_{1,3}^+, 
% s_{z_o,4}=s_{2,1}^+$.  
which are nothing but the basis given in \cite{shioda-usui}.
By the explicit formula for the hight pairing, the Gram matrix with respect
to this basis is the first matrix in Example~\ref{eg:2-nodal}.
We now apply our observation in Section~\ref{sub:bisection}.
 
 Consider
 \begin{align*}
 P = [2]P_0 &=(x_1(t), y_1(t))\\
 &=\left(\frac{1}{144}t^2+\frac{1231}{72}t-\frac{5143775}{144}, -\frac{1}{1728}t^3-\frac{2335}{576}t^2+\frac{13493375}{576}t-\frac{29962489375}{1728}\right)
 \end{align*}
  and $r_{a, 1}(t) = -t/12 + a$. By applying the method 
 given in section~\ref{sub:bisection}, then we have
 \begin{eqnarray*}
 && F(t, x, 1) - \{l_1(t, x)\}^2 \\
 & = &
  -\frac 1{2985984} \left( {t}^{2}+2462\,t-144\,x-5143775 \right) \\
  & & \left( 144\,{a}^{2}{t}^{2}+354528\,{a}^{2}t-20736\,x{a}^{2}+109032\,a{t}^{2}+3456\,xta-740703600\,{a}^{2} 
  \right .\\
  && \left .
\mbox{}-848072400\,ta-86856575\,{t}^{2}-1677600\,tx+20736\,{x}^{2}+719099745000\,a \right .\\
&& \left .+328131278750\,t
\mbox{}+4886010000\,x-174531500609375 \right), 
 \end{eqnarray*}
 where $l_1(t, x) = r_{a, 1}(t)(x - x_1(t)) + y_1(t)$.
 Hence $C_{a,1} = C(r_{a, 1}(t), [2]P_0)$ is a plane curve of degree $2$ given as the zero locus of the second factor.  
 
Next we take 
\begin{align*}
P = P_1\dot{+}P_2 &= (x_2(t), y_2(t))\\
&= \left(\frac{1}{36}t^2+\frac{435}{2}t-\frac{921375}{4},-\frac{1}{216}t^3-\frac{41625}{8}t+\frac{373156875}{8}\right)
\end{align*}
and $r_{b, 2}(t) = -t/6 + b$. By applying the method in 
section~\ref{sub:bisection} again, 
we have
\begin{eqnarray*}
&& F(t, x, 1)^2  - l_2(t, x)^2 \\
& = & -\frac 1{5184}(t^2+7830t-36x-8292375) \\
 & & \left (4b^2t^2+31320b^2t-144b^2x+3732bt^2+48btx-33169500b^2+12555000bt  +683865t^2\right . \\
 && \left .+17208tx+144x^2-13433647500b+1258071750t+5904900x-1360156809375 \right ),
\end{eqnarray*}
where $l_2(t, x) = r_{b, 2}(t)(x - x_2(t)) + y_2(t)$. Hece $C_{b, 2} = C(r_{b, 2}(t), P_1\dot{+}P_2)$ is a 
plane curve of degree $2$ given as the zero locus of the second factor. 

By straightforward computation with computer similar to that in \cite[Lemma 5]{bannai-tokunaga}, one can check that the three condition in the beginning of this section are satisfied. Therefore, we have a Zariski $N+1$-plet.

 \subsection{Case II: $\mcQ$ is an irreducible quartic with one tacnode}

Let $F(T, X, Z)$ be a homogeneous polynomial 
 \[
X^3Z + (25T + 9Z)X^2Z + (144T^2Z + T^3)X + 16T^4 
 \]
and let $\mcQ$ be a quartic given by $F = 0$. $\mcQ$ is an irreducible quartic with one tacnode at $x_o=[0:0:1]$. Let $z_o=[0:1:0]$.
%
%  Let $S$ be the rational elliptic surface whose generic fiber is given by the Weierstrass equation 
%  \[
%  y^2 = F(t, x, 1) = x^3 + (25t + 9)x^2 + (144t^2 + t^3)x + 16t^4
%  \]
 The associated rational elliptic surface $S_{\mcQ,z_o}$ was given and studied in \cite[p. 210]{shioda-usui}.
% it corresponds to considering the rational elliptic surface associated to the quartic $\mathcal{Q}$
%  and $z_o=[0:1:0]$. 
 According to \cite{shioda-usui},
 the rational elliptic surface $S_{\mcQ,z_o}$ has a singular fibers of type $\I_4$ and one singular fiber of type $\III$. 
 The Mordell-Weil lattice of $S$ is  $\MW(S)\cong A_1^\ast\oplus A_3^\ast$ and the narrow Mordell-Weil lattice is  
 $\MW(S)^0\cong A_1\oplus A_3$. 

 Let $L_{x_o}, L_{1}, L_{2}, L_{3}, L_{4}$ be lines defined by
 \begin{align*} 
 L_{x_o}:\, X=0,
\quad  L_{1}:  \,16T+X=0,
\quad   L_{2}: \,15T+X=0,\\
 L_{3}:  \,7T+X=0,
\quad L_{4}:  \,12T+X=0.
 \end{align*}
 These four lines give rise to the following dp-free  sections,
 \begin{align*}
 s_{z_o, 0} =(0, 4t^2),
\quad s_{z_o, 1}=(-16t, -48t),
\quad  s_{z_o,2}=(-15t, -t(t+ 45)),\\
  s_{z_o,3}=(-7t, -3t(t-7)),
\quad s_{z_o,4}=(-12t, 2t(t+ 18)),
 \end{align*} 
 which are nothing but the basis given in \cite{shioda-usui}. We denote the elements corresponding to $s_{z_o, i}$ $(i = 0, \ldots, 3)$ by $P_i$ $(i = 0, \ldots, 3)$, 
 respectively.  
Note that the Gram matrix for $P_0, P_1, P_2, P_3$ is the first one in Example~\ref{eg:tacnode}.
 We have
 \begin{align*}
 [2]P_0 &= (x_1(t), y_1(t)) = \left (\frac 1{64}t^2 - \frac {41}2t + 315, 
 -\frac{55}{512}t^2 + \frac{2637}{8}t - 5670\right ) \\
 P_1\dot{-}P_2 & = (x_2(t), y_2(t)) = \left (t^2 + 192t + 8640, - t^3 -301t^2 - 27936t - 803520 \right )
 \end{align*}

 Take $P = [2]P_0 = (x_1(t), y_1(t))$ and $r_{a,1}(t) = -t/8 + a$. By applying the method in 
Section~\ref{sub:bisection}, we have
\begin{eqnarray*}
&& F(t, x, 1) - \{l_1(t, x)\}^2 \\
& = & - \frac 1{4096}(t^2-1312t-64x+20160)\\
 & & \left (a^2t^2-1312a^2t-64a^2x+548at^2+16atx+20160a^2-47232at+9540t^2 \right .\\
 && \left . +288tx+64x^2+725760a-425088t+20736x+6531840\right ),
 \end{eqnarray*}
  where $l_1(t, x) = r_{a, 1}(t)(x - x_1(t)) + y_1(t)$.
 Hence $C_{a,1} = C(r_{a, 1}(t), [2]P_0)$ is a plane curve of degree $2$ given as the zero locus of the second factor.

 We next take $P = P_1\dot{-}P_2 =(x_2(t), y_2(t))$ and $r_{b, 2}(t) = -t + b$. By Section~\ref{sub:bisection},
 we have 
 \begin{eqnarray*}
 && F(t, x, 1) - \{l(t, x)\}^2 \\
 &= & -(t^2+192t-x+8640)(b^2t^2+192b^2t-b^2x+218bt^2+2btx+8640b^2+38592bt \\
 &&  +11865t^2+217tx+x^2+1607040b+1928448t+8649x+74727360),
 \end{eqnarray*}
 where $l_2(t, x) = r_{b, 2}(t)(x - x_2(t)) + y_2(t)$.
 Hence $C_{b,2} = C(r_{b, 2}(t), P_1\dot{-}P_2)$ is a plane curve of degree $2$ given as the zero locus of the second factor. 
 
By straightforward computation with computer similar to that in \cite[Lemma 5]{bannai-tokunaga}, one can check that the three condition given in the beginning of this section are satisfied.

%Proof of Theorem0.6

\section{Proof of Theorem~\ref{thm:z-5-plet}}
In this section, we prove Theorem~\ref{thm:z-5-plet} by combining the result of Theorem \ref{thm:z-k-plet2} with a known Zariski triple constructed in \cite{bannai16}.

Let $\mcQ$ be the two nodal quartic given in Section \ref{sec:proof_thm:z-k-plet2}. Assume that  there exist $6$ bisections  $D_1, \ldots, D_6$ of $\varphi_{\mcQ, z_o} :
S_{\mcQ, z_o} \to \PP^1$ as follows:

\begin{enumerate}

\item[(i)] $s_{z_o}(D_i) = [2]s_{z_o, i}$ $(i = 1, 2)$, $s_{z_o}(D_i) \in D_4^*$ $(i = 3, 4, 5, 6)$.

\item[(ii)] $\mcC_i = q\circ q_{z_o}\circ f_{\mcQ, z_o}(D_i)$ $(i = 1, \ldots, 5)$ are contact conics
 to $\mcQ$.

\item[(iii)] $\mcC_i$ is tangent to $\mcQ$ at $4$ distinct points for all $i$ and $\mcC_i$ intersects
$\mcC_j$ transversely if $i \neq j$.

\item[(iv)] No three of $\mcC_i$'s meet at one point.

\item[(v)] $\mcC_3, \mcC_4, \mcC_5, \mcC_6$ and $\mcQ$ are the components of  the 
Zariski triple given in  \cite[Theorem~1.4]{bannai16}.

\end{enumerate}

Put
\begin{eqnarray*}
\underline{\mcC}^{(1)} = \mcC_1 + \mcC_2,
\quad\underline{\mcC}^{(2)}  =  \mcC_1 + \mcC_3,
\quad\underline{\mcC}^{(3)}  =  \mcC_3 + \mcC_4 \\
\underline{\mcC}^{(4)}  =  \mcC_3 + \mcC_5,
\quad\underline{\mcC}^{(5)}  =  \mcC_3 + \mcC_6.
\end{eqnarray*}

By \cite[Theorem~1.4]{bannai16}, $(\mcQ + \underline{\mcC}^{(3)}, \mcQ + \underline{\mcC}^{(4)}, \mcQ + \underline{\mcC}^{(5)})$
is a Zariski triple. Note that this Zariski triple was distinguished by considering the \lq\lq splitting type" defined for $(\mcC_i,\mcC_j;\mcQ)$ which takes values $(0,4), (1,3), (2,2)$ in this case. See \cite{bannai16} for details about splitting types.

Next, we consider $\sharp \left ((\Phi^1_{\mcQ, \underline{\mcC}^{(i)}})^{-1}(1)\right )$. We have
\[
\sharp \left ((\Phi^1_{\mcQ, \underline{\mcC}^{(i)}})^{-1}(1)\right )=
\left \{ \begin{array}{cc}
       2 & i = 1 \\
       1 & i = 2 \\
       0 & i = 3, 4, 5
        \end{array} \right.
\]

By combining the two invariants, we have the following table where the top row indicates the value of $\sharp \left ((\Phi^1_{\mcQ, \underline{\mcC}^{(i)}})^{-1}(1)\right )$ and the left most column indicates the splitting type.
\begin{center}
\begin{tabular}{c|ccc}
& 0 & 1 & 2 \\
\hline
$(0,4)$ & $\mcQ+\mcC_3+\mcC_4$ & -- & $\mcQ+\mcC_1+\mcC_2$\\
$(1,3)$ &$\mcQ+\mcC_3+\mcC_5$&--&--\\
$(2,2)$ &$\mcQ+\mcC_3+\mcC_6$& $\mcQ+\mcC_1+\mcC_3$&--
\end{tabular} 
\end{center}
Now it is immediate that $(\mcQ+\underline{\mcC}^{(1)}, \mcQ+\underline{\mcC}^{(2)},\mcQ + \underline{\mcC}^{(3)}, \mcQ + \underline{\mcC}^{(4)}, \mcQ + \underline{\mcC}^{(5)})$ forms a
Zariski $5$-plet. Hence under the assumption that $\mcC_1, \ldots, \mcC_6$ satisfying (i),...,(v) exists, Theorem \ref{thm:z-5-plet} is true.

Finally,we show that $\mcQ$ and $\mcC_i$ ($i = 1, \ldots, 6$) as above exist by providing explicit equations.
Let $\mcC_1,\ldots,\mcC_6\subset\PP^2$ be conics given by $\mcC_1=C(-\frac{1}{12}t,[2]s_0)$, $\mcC_2=C(-\frac{1}{12}t+1,[2]s_0)$, $\mcC_3=C(\frac{1}{20}t,\dot{-}s_2\dot{+}s_3)$, $\mcC_4=C(\frac{1}{20}t+1,\dot{-}s_2\dot{+}s_3)$, $\mcC_5=C(-\frac{1}{24}t,\dot{-}s_1\dot{+}[2]s_2\dot{-}s_3\dot{-}s_4)$, $\mcC_6=C(\frac{1}{6}t,[2]s_1\dot{-}s_2)$. The affine parts of $\mcC_1,\mcC_2$ are given by the following equations:
\begin{align*}
\mcC_1: & {\frac {174531500609375}{20736}}-{\frac {164065639375\,t}{10368}}+{
\frac {86856575\,{t}^{2}}{20736}}-{\frac {33930625\,x}{144}}+{\frac {
5825\,tx}{72}}-{x}^{2}=0\\
\mcC_2: & {\frac {173813141567975}{20736}}-{\frac {163641780439\,t}{10368}}+{
\frac {86747399\,{t}^{2}}{20736}}-{\frac {33930481\,x}{144}}+{\frac {
5813\,tx}{72}}-{x}^{2}=0
%\\
%C_3: & {\frac {11289241}{6400}}\,{t}^{2}-{\frac {
%5829}{40}}\,tx+x^2\\
%&\hspace{2cm}-{\frac {
%914968125}{128}}\,t+{\frac {
%4730609}{16}}\,x+{\frac {1833568580025}{256}}=0\\
%C_4: & {\frac {454825}{256}}\,{t}^{2}-{\frac {
%1165}{8}}\,tx+x^2\\
%&\hspace{2cm}-{\frac {
%914968125}{128}}\,t+{\frac {4730625}{16}}\,x+{\frac {1840331390625}{256}}\\
%C_5: & {\frac {74102665}{4096}}\,{t}^{2}-{\frac {5821}{32}}\,tx+x^2\\
%&\hspace{2cm}-{\frac {
%265613914125}{2048}}\,t+{\frac {40259025}{64}}\,x+{\frac {921634762190625}{4096}}=0\\
%C_6: & {\frac {62785}{16}}\,{t}^{2}-{\frac {361}{2}}\,tx+x^2\\
%&\hspace{2cm}-{\frac {106713693}
%{8}}\,t+{\frac {1225449}{4}}\,x+{\frac {171622907001}{16}}=0\\
\end{align*}
As for the explicit equations of $\mcC_3,\ldots, \mcC_6$, we refer to \cite{bannai16}. It can be easily checked that $\mcC_1,\ldots,\mcC_6$ satisfy the conditions in Section \ref{sec:proofthm:q2}, hence Theorem \ref{thm:z-5-plet} is completely proved.

\begin{rem}
As a final remark, we note that a Zariski 5-plet can be constructed in a similar way in the case where $\mcQ$ is a quartic with one tacnode.
\end{rem}

%\end{document}

%\input example.tex

%\input section4-rev.tex
%
%\input example2.tex

\noindent Shinzo BANNAI\\
Department of Natural Sciences\\
National Institute of Technology, Ibaraki College\\
866 Nakane, Hitachinaka-shi, Ibaraki-Ken 312-8508 JAPAN \\
{\tt sbannai@ge.ibaraki-ct.ac.jp}\\

\noindent Hiro-o TOKUNAGA\\
Department of Mathematics and Information Sciences\\
Graduate School of Science and Engineering,\\
Tokyo Metropolitan University\\
1-1 Minami-Ohsawa, Hachiohji 192-0397 JAPAN \\
{\tt tokunaga@tmu.ac.jp}

\vspace{0.5cm}
      
%\noindent Masayuki KAWASHIMA\\
%      Department of Mathematics,\\
%         Tokyo University of Science,\\
%         1-3 Kagurazaka, Shinjuku-ku, 
%         Tokyo 162-8601 JAPAN\\
%{\tt kawashima@ma.kagu.tus.ac.jp}
% \input alexander1.tex     

\begin{thebibliography}{99}
\bibitem{Artal}
E.~Artal~Bartolo:
\newblock \emph{Sur les couples des {Zariski}},
\newblock {\rm J.\ Algebraic Geometry}, {\bf 3} (1994) no. {\bf 2}, 223--247

   \bibitem{act} E.~Artal Bartolo, J.-I.~Cogolludo and H.~Tokunaga:
   \emph{A survey on Zariski pairs}, Adv. Stud. Pure Math., \textbf{50} (2008), 1-100.
   
   
  
  
  %
  %
  %
  
  \bibitem{bannai16}  S.~Bannai: \emph{A note on splitting curves of plane quartics and multi-sections of rational elliptic surfaces}, Topology and its Applications {\bf 202} (2016), 428-439.

   
   \bibitem{bgst} S.~Bannai, B.Guerville-Ball\'{e}, T.~Shirane and H.~Tokunaga:
  \emph{On the topology of arrangements of a cubic and its inflectional tangents}, arXiv:1607.07618.
  
  \bibitem{bannai-shirane16} S.~Bannai and T.~Shirane: \emph{Nodal curves with a contact-conic and 
  Zariski pairs}, arXiv:1608.03760.
  
 \bibitem{bannai-tokunaga} S.~Bannai and H.~Tokunaga:  \emph{Geometry of  bisections of elliptic surfaces  and 
 Zariski $N$-plets for conic arrangements}, Geom.  Dedicata
 {\bf  178} (2015),  219-237,  DOI 10.1007/s10711-015-0054-z. 


%\bibitem{bannai-kawashima-tokunaga} Bannai, S. ,  Kawashima, M. and Tokunaga, H.: \emph{On the topology of the complements of reducible plane curves via Galois covers}, arXiv:1304.0536
 
% \bibitem{bpv} W.~Barth, K.~Hulek, C.A.M.~Peters and A. Van de Ven: Compact complex surfaces, 
% Ergebnisse der Mathematik und ihrer Grenzgebiete {\bf 4} 2nd Enlarged Edition, Springer-Verlag (2004).
% 
%\bibitem{CLO} D.~Cox, J.~Little, D.~O'Shea: Ideals, varieties and algorithms. An introduction to computational algebraic geometry and commutative algebra. 3rd edition, Springer-Verlag (2007).

%\bibitem{Degtyarev-alex}
%A.~I. Degtyarev.
%\newblock \emph{Alexander polynomial of a curve of degree six},
%\newblock {\rm J. Knot Theory Ramifications}, {\bf 3} (1994), 439--454.

%\bibitem{Eisenbud-Neumann} D.~Eisenbud and W. Neumann: Three-dimensional link theory and invariants of plane
%curve singularities, Ann. of Math. Stud. {\bf 110}, Princeton Univ. Press (1985).

%\bibitem{Esnault}
%H.~Esnault:
%\newblock \emph{Fibre de {M}ilnor d'un c\^one sur une courbe plane singuli\`ere},
%\newblock {\rm Invent. Math.}, {\bf 68} (1982) no. {\bf 3}, 477--496.

\bibitem{fulton} W.~Fulton: Intersection Theory, Springer-Verlag (1984).

 \bibitem{horikawa} E.~Horikawa: \emph{ On deformations of quintic surfaces},
\rm Invent. Math. {\bf 31} (1975), \rm $43 - 85$.

%\bibitem{Kawashima2}
%M.~Kawashima and M.~Oka:
%\newblock \emph{On {A}lexander polynomials of certain {$(2,5)$} torus curves.}
%\newblock {\rm J. Math. Soc. Japan}, {\bf 62} (2010) no. {\bf 1}, 213--238.


\bibitem{kodaira} K.~Kodaira: \emph{On compact analytic surfaces II-III}, Ann. of Math. \textbf{77}
(1963), 563-626, \textbf{78} (1963), 1-40.
%\bibitem{MR85h:14017}
%A.~Libgober:
%\newblock \emph{Alexander polynomial of plane algebraic curves and cyclic multiple planes.}
%\newblock Duke Math. J. \textbf{49} (1982), no. 4, 833--851
%\bibitem{Loeser-Vaquie}
%F.~Loeser and M.~Vaqui{\'e}.
%\newblock \emph{Le polyn\^ome d'{A}lexander d'une courbe plane projective},
%\newblock {\rm Topology}, {\bf 29} (1990) no. {\bf 2}, 163--173.
  
%\bibitem{miranda} R.~Miranda: \emph{The moduli of Weierstrass fibrations over ${\mathbb P}^1$}, 
% Math. Ann. \textbf{255}(1981), 379-394.
 
\bibitem{miranda-basic} R.~Miranda: \emph{Basic theory of elliptic surfaces}, Dottorato di Ricerca in Matematica, ETS Editrice, Pisa, 1989.

% 
\bibitem{miranda-persson} R.~Miranda and U.~Persson: \emph{
On extremal rational elliptic surfaces}, Math. Z. \textbf{193} (1986), 537-558.

%
%\bibitem{namba-tsuchi}
%M.~Namba and H.~Tsuchihashi, \emph{On the fundamental groups of {G}alois
%  covering spaces of the projective plane}, Geom. Dedicata~\textbf{104} (2004),
%  97--117.
%
\bibitem{oguiso-shioda}  K.~Oguiso and T.~Shioda: \emph{The Mordell-Weil lattice of a Rational
Elliptic surface}, Comment. Math. Univ. St. Pauli \textbf{40} (1991), 83-99.

%\bibitem{Okabook}
%M.~Oka.
%\newblock \emph{Non-degenerate complete intersection singularity}.
%\newblock Hermann, Paris, 1997.
%
%\bibitem{OkaAtlas}
%M.~Oka.
%\newblock \emph{Alexander polynomial of sextics.}
%\newblock {\rm J. Knot Theory Ramifications}, {\bf 12} (2003) no. {\bf 5}, 619--636.
%
%\bibitem{OkaSurvey}
%M.~Oka.
%\newblock \emph{A survey on {Alexander} polynomials of plane curves.}
%\newblock {\rm Singularit\'es Franco-Japonaise, S\'eminaire et congr\`es},
%  {\bf 10} Soc. Math. France, Paris, 2005.
%%
%\bibitem{AEC} J.H.~Silverman: \emph{The Arithmetic of Elliptic Curves}, Graduate Texts in Mathematics, {\bf 106}
%Springer-Verlag, 1985.
%
%\bibitem{silverman} J.H.~Silverman: \emph{Advanced Topics in the Arithmetic of Elliptic Curves}, Graduate Texts in Mathematics, {\bf151} Springer-Verlag, 1994
%

 \bibitem{shioda90} T.~Shioda: \emph{On the Mordell-Weil lattices}, \rm Comment. Math. Univ. St. Pauli
\textbf{39} (1990), 211-240.  


\bibitem{shioda92} T.~Shioda: \emph{Existence of a rational elliptic surface with a given Mordell-Weil lattice},  
\rm Proc. Japan Acad. Ser. A Math. Sci. {\bf 68} (1992), no. {\bf 9}, 251--255. 

\bibitem{shioda93} T.~Shioda: \emph{Plane Quartics and Mordell-Weil Lattices of Type $E_7$},
\rm Comment. Math. Univ. St. Pauli
\textbf{42} (1993), 61--79.  


 \bibitem{shioda-usui} T.~Shioda and H.~Usui: \emph{Fundamental invariants of Weyl groups and
 excellent families of elliptic curves}, Comment. Math. Univ. St. Pauli \textbf{41} (1992),
 169-217.
% 

\bibitem{shirane16} T.~Shirane: \emph{A note on splitting numbers for
Galois covers and $\pi_1$-equivalent Zariski $k$-plets}, Proc. AMS., DOI 10.1090/proc/13298



%\bibitem{tokunaga94} H.~Tokunaga:  \emph{On dihedral Galois coverings}, \rm Canadian J. of
%Math. {\bf 46} \rm (1994),1299 - 1317.
%
%\bibitem{tokunaga97} H.~Tokunaga: \emph{Dihedral coverings branched along maximizing sextics}, \rm Math. Ann. {\bf 308}  (1997)
%633-648.
%%%
%\bibitem{tokunaga98}  H.~Tokunaga: \emph{Some examples of Zariski pairs arising from certain 
%elliptic K3 surfaces I}, \rm
%Math. Z. {\bf 227} (1998), 465-477, \emph{Some examples of Zariski pairs arising from certain 
%elliptic K3 surfaces II: Degtyarev's conjecture}, 
% Math. Z. {\bf 230} (1999), 389-400
%
%
%\bibitem{tokunaga04} H.~Tokunaga: \emph{Dihedral covers and an elementary arithmetic on elliptic surfaces}, 
%J. Math. Kyoto Univ. \textbf{44}(2004), 55-270.

\bibitem{tokunaga10} H.~Tokunaga: \emph{Geometry of irreducible plane quartics  and  their quadratic residue conics}, 
 J. of Singularities(electric), \textbf{2} (2010), 170-190.

\bibitem{tokunaga12} H.~Tokunaga: \emph{Some sections on rational elliptic surfaces and certain special
conic-quartic configurations}, Kodai Math. J.\textbf{35} (2012), 78-104.


\bibitem{tokunaga14} H.~Tokunaga: \emph{
 Sections of elliptic surfaces 
and
 Zariski pairs for conic-line arrangements via dihedral covers
 }, J. Math. Soc. Japan {\bf 66} (2014), 613-640.

\bibitem{tumenbayar-tokunaga} K.~Tumenbayar and H.~Tokunaga: \emph{Elliptic surfaces and contact conics for a 
$3$-nodal quartic}, to appear in Hokkaido Math. J.

%
%%
%\bibitem{zariski29} O.~Zariski: 
%\emph{On the problem of existence of algebraic functions of two variables possessing a 
%given branch curve}, Amer. J. Math.~\textbf{51} (1929), 305--328.
%%
%%
%\bibitem{zariski37}
%O.~Zariski:  \emph{The topological discriminant group of a {R}iemann surface of genus $p$}, 
%Amer. J. Math.~\textbf{59} (1937), 335--358.
%%
%\bibitem{zariski} O. Zariski: \emph{On the purity of the branch locus of algebraic functions}, 
%Proc. Nat. Acad. USA \textbf {44} (1958), 791-796.
%




\end{thebibliography}
 \end{document}